\documentclass[12pt]{article}

\usepackage[left=2.7cm,bottom=2.7cm,right=2.7cm,top=2.7cm]{geometry}

\usepackage{multirow}
\usepackage{booktabs}
\usepackage{amsmath,amsfonts}
\usepackage{array}
\usepackage[caption=false,font=normalsize,labelfont=sf,textfont=sf]{subfig}
\usepackage{textcomp}
\usepackage{stfloats}
\usepackage{url}
\usepackage{verbatim}
\usepackage{graphicx}
\usepackage{hyperref}
\usepackage{amssymb}
\usepackage{mathtools}
\usepackage{algorithm}
\usepackage{algorithmic}
\usepackage{cases}
\usepackage{bbding}
\usepackage[bottom]{footmisc}

\newtheorem{theorem}{Theorem}
\newtheorem{corollary}{Corollary}

\newtheorem{remark}{Remark}
\newtheorem{assumption}{Assumption}
\newtheorem{proposition}{Proposition}

\usepackage[utf8]{inputenc}
\usepackage{tikz}
\usetikzlibrary{patterns}
\usetikzlibrary{decorations.pathreplacing,calligraphy}
\usetikzlibrary {patterns.meta}
\usetikzlibrary {arrows.meta}
\usetikzlibrary{calc}
\tikzdeclarepattern{
  name=mylines,
  parameters={
      \pgfkeysvalueof{/pgf/pattern keys/size},
      \pgfkeysvalueof{/pgf/pattern keys/angle},
      \pgfkeysvalueof{/pgf/pattern keys/line width},
  },
  bounding box={
    (0,-0.5*\pgfkeysvalueof{/pgf/pattern keys/line width}) and
    (\pgfkeysvalueof{/pgf/pattern keys/size},
0.5*\pgfkeysvalueof{/pgf/pattern keys/line width})},
  tile size={(\pgfkeysvalueof{/pgf/pattern keys/size},
\pgfkeysvalueof{/pgf/pattern keys/size})},
  tile transformation={rotate=\pgfkeysvalueof{/pgf/pattern keys/angle}},
  defaults={
    size/.initial=5pt,
    angle/.initial=45,
    line width/.initial=.4pt,
  },
  code={
      \draw [line width=\pgfkeysvalueof{/pgf/pattern keys/line width}]
        (0,0) -- (\pgfkeysvalueof{/pgf/pattern keys/size},0);
  },
}

\usepackage{amsmath,amsfonts}
\usepackage{array}
\usepackage[caption=false,font=normalsize,labelfont=sf,textfont=sf]{subfig}
\usepackage{textcomp}
\usepackage{stfloats}
\usepackage{url}
\usepackage{verbatim}
\usepackage{graphicx}
\usepackage{mathtools}
\usepackage{algorithm}
\usepackage{algorithmic}
\usepackage{cases}
\usepackage{enumerate}
\usepackage{tabularx}
\usepackage{multirow}
\usepackage{algorithm}
\usepackage{authblk}
\usepackage{algorithmic}
\usepackage{indentfirst}
\usepackage{setspace}
\usepackage{color}
\usepackage[utf8]{inputenc} 
\usepackage[T1]{fontenc}    
\usepackage{hyperref}       
\usepackage{url}            
\usepackage{booktabs}       
\usepackage{amsfonts}       
\usepackage{nicefrac}       
\usepackage{microtype}      
\usepackage{lipsum}		
\usepackage{graphicx}
\usepackage{natbib}
\usepackage{doi}
\usepackage{amssymb}                  
\usepackage{mathrsfs} 
\usepackage[resetlabels]{multibib}
\usepackage[bottom]{footmisc}

\usepackage[utf8]{inputenc}
\usepackage{tikz}

\begin{document}

\title{Rethinking Hard Thresholding Pursuit: Full Adaptation and Sharp Estimation}

\author[1]{Yanhang Zhang}
\author[2]{Zhifan Li}
\author[1]{Shixiang Liu}
\author[3]{Xueqin Wang}
\author[1, 4]{Jianxin Yin}

\affil[1]{\footnotesize School of Statistics, Renmin University of China}
\affil[2]{\footnotesize Beijing Institute of Mathematical Sciences and Applications}
\affil[3]{\footnotesize Department of Statistics and Finance/International Institute of Finance, School of Management, University of Science and Technology of China}
\affil[4]{\footnotesize Center for Applied Statistics and School of Statistics, Renmin University of China}

\date{}
\maketitle \sloppy

\begin{abstract}
  Hard Thresholding Pursuit (HTP) has aroused increasing attention for its robust theoretical guarantees and impressive numerical performance in non-convex optimization.
  In this paper, we introduce a novel tuning-free procedure, named Full-Adaptive HTP (FAHTP), that simultaneously adapts to both the unknown sparsity and signal strength of the underlying model. 
 We provide an in-depth analysis of the iterative thresholding dynamics of FAHTP, offering refined theoretical insights.  In specific, under the beta-min condition $\min_{i \in S^*}|{\boldsymbol{\beta}}^*_i| \ge C\sigma (\log p/n)^{1/2}$, we show that the FAHTP achieves oracle estimation rate $\sigma (s^*/n)^{1/2}$, highlighting its theoretical superiority over convex competitors such as LASSO and SLOPE, and recovers the true support set exactly. More importantly, even without the beta-min condition, our method achieves a tighter error bound than the classical minimax rate with high probability.
 The comprehensive numerical experiments substantiate our theoretical findings, underscoring the effectiveness and robustness of the proposed FAHTP.
\end{abstract}


\section{Introduction} 
Consider a linear model 
\begin{equation*}
 \boldsymbol{y} =  \boldsymbol{X}{\boldsymbol{\beta}}^* + \boldsymbol{\xi},
\end{equation*}
where $\boldsymbol{y}\in \mathbb{R}^n$, $\boldsymbol{X} \in \mathbb{R}^{n\times p}$ and ${\boldsymbol{\beta}}^* \in \mathbb{R}^p$ are the response vector, the design matrix and the underlying regression coefficient respectively. And $\boldsymbol{\xi} \in \mathbb{R}^n$ is the sub-Gaussian random error with scale parameter $\sigma^2$. 
In the high-dimensional learning framework, we focus on the case where $p \gg n$ and the coefficient ${\boldsymbol{\beta}}^*$ is $s^*$-sparse.

One typical learning approach under empirical risk minimization is the best-subset selection problem, formulated as
\begin{equation}\label{bess}
 \min_{{\boldsymbol{\beta}} \in \mathbb{R}^p} \mathcal{L}_n( {\boldsymbol{\beta}}), \quad s.t.\ \| {\boldsymbol{\beta}}\|_0 \leq s,
\end{equation}
where $\mathcal{L}_n( {\boldsymbol{\beta}}) = \| \boldsymbol{y}- \boldsymbol{X} {\boldsymbol{\beta}}\|_2^2/(2n)$ and $s$ is some given positive integer controlling the sparsity level. However, problem \eqref{bess} is an NP-hard problem \citep{natarajan1995sparse}, indicating it lacks a polynomial-time solution.

Among the numerous solutions that aim at approximately solving the $\ell_0$-sparse problem \eqref{bess}, the Iterative Hard Thresholding (IHT, \cite{blumensath2009iterative}) algorithm stands out as an extensively studied iterative greedy selection approach.
It keeps the $s$ largest values of the gradient descent in each step via a hard thresholding procedure, yielding an estimator with $s$ nonzero components.
To enhance empirical performance, \cite{foucart2011hard} introduced a debiasing step at each iteration, employing an ordinary least-squares (OLS) estimator—a method they named Hard Thresholding Pursuit (HTP).
The extensive literature on HTP and its variants includes both statistics and optimization \citep{ zhang2018,huang2018constructive, Zhu202014241}. In practice, HTP-style methods have found their wide applications in deep neural networks pruning \citep{benbaki2023fast}, genetic data analysis \citep{zhang2022, zhang2024} and Alzheimer’s Disease Neuroimaging Initiative study \citep{hao2021optimal}.
In theory, \cite{jain2014iterative} relaxed the bounded restriction on the RIP-type condition for $s\gg s^*$, and showed the minimax optimality for parameter estimation of IHT up to a logarithm factor.
\cite{zhang2018} analyzed the parameter estimation and support recovery of HTP for both $s=s^*$ and $s\gg s^*$.
Additionally, \cite{huang2018constructive} established the $\ell_2$ and $\ell_\infty$ estimation error bounds for $s\geq s^*$. 

\subsection{Main results}
Our work focuses on the estimation and minimax adaptation properties of HTP. In contrast to techniques based on complexity theory and empirical processes used in convex estimators like LASSO \citep{bellec2018SLOPE}, or non-convex algorithms such as SCAD and MCP \citep{loh2013regularized,fan2014}, our approach offers a detailed analysis of the iterative thresholding dynamics, providing a deeper understanding of the intrinsic properties of the hard thresholding estimator. By relying solely on the RIP assumption \eqref{rip}, we demonstrate that HTP algorithm has the following key properties:
\begin{proposition}\label{pro1}
    We divide the true support set $S^*$ into 
\begin{equation*}
   S_1^*\coloneqq\{i \in S^*: |{\boldsymbol{\beta}}^*_i| \geq C_1 \sigma(\log p/n)^{1/2}\},\quad S_2^* \coloneqq S^* \setminus S_1^*.
\end{equation*}
Assume
    $|S_1^*|\log p/{n} \to 0, ~ {|S_1^*|}/{p} \to 0$ and $|S_1^*| \to \infty$.
 Let $s$ be the input parameter in Algorithm \ref{alg:iht1}. Given $s = |S_1^*|$, with probability tending to 1, the estimator of HTP has
 \begin{equation}\label{adaptive}
\|\hat{\boldsymbol{\beta}} - \boldsymbol{\beta}^*\|_2 \leq C_2\sigma\left(\frac{|S_1^*|}{n}\right)^{1/2} + C_3\|\boldsymbol{\beta}^*_{S_2^*}\|_2.
\end{equation}
\end{proposition}
Proposition \ref{pro1} shows that HTP can achieve a tighter estimation upper bound better than the classical minimax rate $\sigma(s^*\log p/n)^{1/2}$. Specifically, HTP successfully detects the large signal set $S_{1}^*$ and provides unbiased estimates for these signals, while the hard thresholding operator simultaneously shrinks the small signals in $S_2^*$. This combination of advantages allows HTP to surpass the classical minimax rate and derive a sharper error bound \eqref{adaptive}. 

Based on the estimation results in Proposition \ref{pro1}, we propose a novel tuning-free strategy (FAHTP, Algorithm \ref{alg:iht2}) for selecting the sparsity parameter $s$.
\begin{proposition}\label{pro2}
    Assume that $\min_{i \in S^*}|{\boldsymbol{\beta}}^*_i| \geq C_1 \sigma(\log p/n)^{1/2}$,
and ${s^*\log p}/{n} \to 0,~{s^*}/{p} \to 0,~ s^* \to \infty
    $.
Then, by the proposed FAHTP, with probability tending to 1, we have
 \begin{equation*}\label{oracle}
     \|\hat{\boldsymbol{\beta}} - \boldsymbol{\beta}^*\|_2 \le C_2\sigma\left(\frac{s^*}{n}\right)^{1/2}
 \end{equation*}
 and $\boldsymbol{\hat\beta}$ can recover the true support set $S^*$ exactly.
\end{proposition}
Proposition \ref{pro2} demonstrates that FAHTP achieves exact support recovery under the optimal beta-min condition, as established in \cite{butucea2018variable}. These results indicate that FAHTP successfully inherits the strengths of the hard thresholding operator in the context of the Gaussian sequence model \citep{butucea2018variable, john17gaus} under the RIP assumption.

\subsection{Related literature and comparison}

 Over the past few decades, convex relaxation algorithms to the problem (\ref{bess}) such as LASSO \citep{T1996} and SLOPE \citep{bogdan2015SLOPE} have been widely studied. While these convex algorithms can achieve the minimax optimal rate $\sigma\{s^*\log (p/s^*)/n\}^{1/2}$ \citep{bickel2009simultaneous, bellec2018SLOPE}, they may not be admissible in certain circumstances. 
In the Gaussian sequential model, \cite{ndaoud2019interplay} uncovered a phase transition phenomenon in the minimax rates under the beta-min condition $\min_{i \in S^*}|{\boldsymbol{\beta}}^*_i| \ge C\sigma \{\log (p/s^*)/n\}^{1/2}$. 
They demonstrated that, under this condition, the minimax rate improves to the oracle estimation rate  $\sigma(s^*/n)^{1/2}$, which is independent of the dimension $p$. This finding suggests that an appropriate beta-min condition can lead to sharper estimation error bounds. However, \cite{bellec2018noise} showed that convex algorithms can only achieve the original minimax optimal rate when the beta-min condition holds. 

The above observations lead to an important question: when and why does the HTP algorithm outperform the convex counterparts? To be specific, we aim to investigate when the HTP can exceed the minimax rate and achieve the oracle estimation rate 
$\sigma(s^*/n)^{1/2}$.

 \cite{Zhu202014241, yuan2022stability} investigated the theoretical guarantees of the HTP-style method under a stronger beta-min condition. Specifically, \cite{Zhu202014241} proposed a tuning-free procedure using a splicing technique and demonstrated that it could recover the support set $S^*$ exactly under a stronger condition $\min_{i \in S^*}|{\boldsymbol{\beta}}^*_i| \ge C\sigma \{s^*\log p\log\log n/n\}^{1/2}$.
 \cite{yuan2022stability} established the fast rates for the empirical risk of IHT relying on condition $\min_{i \in S^*}|{\boldsymbol{\beta}}^*_i| \ge C\sigma \{s^*\log p/n\}^{1/2}$.
 Recent work \citep{roy2022high,zhu2024sure} proved the variable selection and FDR-control behavior of IHT-style methods based on condition $\min_{i \in S^*}|{\boldsymbol{\beta}}^*_i| \ge C\sigma (\log p/n)^{1/2}$. These studies also assume $s^* = O(\log p)$, a condition not required in our analysis. Notably, the aforementioned works did not establish the oracle estimation rate for HTP-style methods.

Importantly, the minimax adaptive procedure for unknown $s^*$ of HTP-style methods still remains unexplored. \cite{birge2001gaussian} conducted the Gaussian model selection based on complexity penalization with a known noise level $\sigma$. Building on this approach, \cite{verzelen2012minimax} proposed a novel information criterion to achieve minimax adaptivity for both unknown $s^*$ and $\sigma$. However, this procedure relies on an exhaustive search with subset sizes less than $[n/4]$, making it inefficient for practical implementation. Therefore, developing a minimax adaptive procedure using efficient methods like HTP is a significant direction for high-dimensional model selection.
 
Thus, the main drawbacks in the existing analysis for HTP are listed as follows:
\begin{itemize}
    \item Firstly, the signal condition required in \cite{Zhu202014241,yuan2022stability} is relatively stronger than the condition(see in Theorem \ref{thm:scaled}) listed in this paper, which is the optimal beta-min condition proved by \cite{ndaoud2019interplay}.
    \item Secondly, whether the classical HTP-style algorithm can achieve the oracle estimation rate $\sigma (s/n)^{1/2}$ remains unclear, indicating further analysis is needed.
    \item Lastly, previous studies about HTP (e.g.,\cite{huang2018constructive,zhang2018}) does not imply a minimax adaptive procedure for sparsity $s^*$.
\end{itemize}

Some of the works mentioned below theoretically inspire us. \cite{butucea2018variable} analyzed variable selection from the non-asymptotic perspective with hamming loss, and this can be linked to estimation risk \citep{ndaoud2019interplay}. 
In recent work, \cite{ndaoud2020scaled} introduced a novel IHT-style procedure, extending the theory from \cite{ndaoud2019interplay, butucea2018variable} to linear regression. This procedure involves dynamically updating the threshold geometrically at each step until it reaches a universal statistical threshold. Our procedure, choosing the sparsity $s$ exactly instead of dynamically updating the threshold, offers one main practical 
advantage as follows. Because the model size is discrete while the tuning parameter is continuous, the tuning for model size $s$ is more intuitive in practice and is computationally more efficient compared to the M-estimators.  

\subsection{Main contributions}
We summarize the main contributions of our paper as the following two points:
\begin{itemize}
    \item Under the beta-min condition $\min_{i \in S^*}|{\boldsymbol{\beta}}^*_i| \ge C\sigma \{\log (p/s^*)/n\}^{1/2}$, our study reveals that the HTP algorithm surpasses the minimax rate $\sigma \{s^*\log (p/s^*)/n\}^{1/2}$ and attains the oracle estimation rate $\sigma (s^*/n)^{1/2}$. 
    Moreover, we characterize the dynamic signal detection properties of the HTP algorithm. Because the hard thresholding performs signal detection through iterations, the algorithm can shrink signals according to magnitudes, which leads to the adaptation to signals as well as the oracle estimation rate.
   
    \item  We propose a novel tuning-free procedure (FAHTP) to adapt to the unknown sparsity $s^*$. In the first step, the information criterion \citep{verzelen2012minimax} ensures that HTP achieves minimax optimality with unknown $s^*$ and $\sigma$.
    In the second step, under beta-min condition $\min_{i \in S^*}|{\boldsymbol{\beta}}^*_i| \ge C\sigma (\log p/n)^{1/2}$, FAHTP can further attain the oracle estimation rate and recover the support set $S^*$ exactly with high probability. Moreover, when the beta-min condition is not satisfied, FAHTP can still derive a tighter error bound than the minimax rate with high probability.

\end{itemize}

Based on the above two points, solely under RIP assumption, our proposed algorithm is minimax adaptive to the signal strength and support size $s^*$, which fully matches the classical minimax rate and oracle rate \citep{ndaoud2019interplay}. 

\subsection{Organization}
The present paper is organized as follows. In Section \ref{sec:algorithm}, we establish the estimation error bounds of HTP when $s\geq s^*$, which matches the known minimax lower bounds and demonstrates the optimality of HTP. Based on these results, we propose a minimax adaptive procedure with unknown $s^*$ and $\sigma$.
In Section \ref{sec:scaled}, under the beta-min condition, we prove that the HTP estimator can attain the oracle estimation rate and can achieve the exact recovery of the true support set when $s = s^*$. In Section \ref{sec:strong_signal}, we extend the analysis to general signals, i.e., without the beta-min conditions, and show that our estimator can obtain tighter error bounds with high probability.
Section \ref{sec:simulation} presents numerical experiments that illustrate our theoretical findings.
Lastly, our conclusions are drawn in Section \ref{sec:conclusion}. Detailed proofs of the theory and additional simulations are available in the appendix. More importantly, we provide a road map of the proof of the main results in the appendix.
\subsection{Notations}
For the given sequences $a_n$ and $b_n$, we say that $a_n = O(b_n)$ (resp. $a_n  = \Omega(b_n)$) when $a_n \le cb_n$ (resp. $a_n \ge c b_n$) for some
positive constant $c$. We write $a_n \asymp b_n$ if $a_n = O(b_n)$ and $a_n  = \Omega(b_n)$ while $a_n = o(b_n)$ corresponds to $a_n/b_n \rightarrow 0$ as $n$ goes to infinity.
Denote $[m]$ as the set $\{1,2,\ldots,m\}$, and $\mathrm{I}(\cdot)$ as the indicator function.
Let $x \vee y $ be the maximum of $x$ and $y$, while $x \wedge y $ is the minimum of $x$ and $y$.
Denote $S^* = \{i: {\boldsymbol{\beta}}^*_i \neq 0\} \subseteq [p]$ as the support set of ${\boldsymbol{\beta}}^*$.
For any set $S$ with cardinality $|S|$,
let ${\boldsymbol{\beta}}_S=({\boldsymbol{\beta}}_j, j \in S) \in \mathbb{R}^{|S|}$ and $ \boldsymbol{X}_S = ( \boldsymbol{X}_j, j \in S) \in \mathbb{R}^{n \times |S|}$, and let $( \boldsymbol{X}^\top  \boldsymbol{X})_{SS} \in \mathbb{R}^{|S|\times |S|}$ be the submatrix of $ \boldsymbol{X}^\top  \boldsymbol{X}$ whose rows and columns are both listed in $S$.
For a vector ${\boldsymbol{\beta}}$, denote $\|{\boldsymbol{\beta}}\|_2$ as its Euclidean norm and $\|{\boldsymbol{\beta}}\|_0$ as the count of the nonzero elements of ${\boldsymbol{\beta}}$. Let $\text{supp}({\boldsymbol{\beta}})$ be the support set of ${\boldsymbol{\beta}}$.
For a matrix $ \boldsymbol{X}$,
denote $\| \boldsymbol{X}\|_2$ as its spectral norm and $\| \boldsymbol{X}\|_F$ as its Frobenius norm. Denote $\mathbb{I}_p$ as the $p\times p$ identity matrix.
For two sets $A$ and $B$, define the symmetric difference of $A$ and $B$ as $A \Delta B=(A \backslash B) \cup(B \backslash A)$.
Let $C, C_0, C_1,\ldots$ denote positive constants whose actual values vary from time to time.
To facilitate computation, we assume $\| \boldsymbol{X}_j\|_2  = n^{1/2}$ for all $j \in [p]$.

Notably, the asymptotic results are considered as $p \rightarrow \infty$ when all other parameters of the problem, i.e., $n, s^*$, depend on $p$ in such a way that $n = n(p) \rightarrow \infty$.
For brevity, the dependence of these parameters on $p$ will be further omitted in the notation. 
\section{Minimax optimality and adaptation of HTP }\label{sec:algorithm} 
In this section, we present the nonasymptotic error bounds for the HTP estimator. We first provide the HTP algorithm in Section \ref{sec:IHT}. Then, the $\ell_2$ error bounds are established in Section \ref{sec:l2}. Based on these results, Section \ref{sec:adaptive} develops a minimax adaptive procedure adaptive to the unknown sparsity and noise level.
\subsection{HTP algorithm}\label{sec:IHT} 
To enforce the sparsity of the solution, HTP applies a hard thresholding operator $\mathcal{T}_s(\cdot)$ to the gradient descent in each iteration, where the operator $\mathcal{T}_s({\boldsymbol{\beta}})$ is a non-linear operator that sets all but the largest (in absolute sense) $s$ elements of ${\boldsymbol{\beta}}$ to zero. Then, it incorporates a debiasing step to each step. 
The HTP algorithm is presented in Algorithm \ref{alg:iht1}.

\begin{algorithm}[h]
  \caption{\label{alg:iht1}\textbf{H}ard \textbf{T}hresholding \textbf{P}ursuit (HTP) algorithm.}
  \begin{algorithmic}[1]
    \REQUIRE $\boldsymbol{X}$, $\boldsymbol{y}$ and $s$.\\
    \STATE Initialize $t=0$, ${\boldsymbol{\beta}}^t = \mathbf 0$ and $S^t = \emptyset$.
    \WHILE {$S^{t+1}\neq S^t, $}
    \STATE ${\tilde{\boldsymbol{\beta}}}^{t+1} = \mathcal{T}_{s}\left({{\boldsymbol{\beta}}}^{t} +   \boldsymbol{X}^{\top}( \boldsymbol{y}- \boldsymbol{X}{{\boldsymbol{\beta}}}^{t})/n\right)$.
    \STATE $S^{t+1} = \text{supp}(\tilde{\boldsymbol{\beta}}^{t+1})$.
    \STATE ${\boldsymbol{\beta}}^{t+1} = \arg\min\{\|\boldsymbol{y}-\boldsymbol{X}{\boldsymbol{\beta}}\|_2,\ \text{supp}({\boldsymbol{\beta}}) = S^{t+1}\}.$
    \STATE $t = t+1$.
    \ENDWHILE
    \ENSURE $\hat {\boldsymbol{\beta}} = {\boldsymbol{\beta}}^{t}$.
  \end{algorithmic}
\end{algorithm}
We can decompose the iterative term into the following parts: 
\begin{align*}
  {\boldsymbol{\beta}}^{t} +  \boldsymbol{X}^{\top} ( \boldsymbol{y}- \boldsymbol{X}{\boldsymbol{\beta}}^{t})/n
  = & {\boldsymbol{\beta}}^* + \underbrace{{\boldsymbol{\Phi}} ({\boldsymbol{\beta}}^*-{\boldsymbol{\beta}}^t)}_{\text {optimization error }} + \underbrace{\boldsymbol{\Xi}}_{\text {statistical error }},    
\end{align*}
where 
${\boldsymbol{\Phi}}\coloneqq  \boldsymbol{X}^{\top} \boldsymbol{X}/n-\mathbb{I}_p\ \text{and}\ \boldsymbol{\Xi}  \coloneqq  \boldsymbol{X}^{\top}\boldsymbol{\xi}/n .$
From the decomposition, we can have a view that the error at each step mainly comes from two parts, that is, the estimation error incurred by ${\boldsymbol{\beta}}^t$ and the statistical error incurred by the white noise $\boldsymbol{\xi}$. In the subsequent section, we discuss how to control these two sources of errors during iterations.
\subsection{\texorpdfstring{$\ell_2$}- error bounds}\label{sec:l2} 
We say that the design matrix $ \boldsymbol{X}$ satisfies Restricted Isometry Property (RIP) with order $s$ and constant $0 < \delta_s < 1$ if for all $ S \subset [p]$ with $|S| \leq s$ and $\forall {u} \neq 0, {u} \in \mathbb{R}^{|S|}$, it holds that
\begin{equation*}
1-\delta_s \leq \frac{\left\|{\boldsymbol{X}}_{S} {u}\right\|_2^2}{n\|{u}\|_2^2} \leq 1+\delta_s.
\end{equation*}
Denote this condition by ${\boldsymbol{X}} \sim \operatorname{RIP}(s, \delta_s)$.
We present assumptions below to establish the theoretical properties of the HTP  algorithm.

\begin{assumption}\label{subgaussian}
The random errors ${ \boldsymbol{\xi}}_1,\ldots,{ \boldsymbol{\xi}}_n$ are i.i.d with mean zero and sub-Gaussian tails, that is, there exists a positive number $\sigma$ such that 
\begin{equation*}
{\rm{pr}}(|{ \boldsymbol{\xi}}_i|>t) \leqslant 2\exp\left(-\frac{t^2}{2\sigma^2} \right), \quad \text{for all } t \geqslant 0.
\end{equation*} 
\end{assumption}
\begin{assumption}\label{rip}
 Design matrix $ \boldsymbol{X}$ satisfies \mbox{RIP}$(3s, \delta_{3s})$. \end{assumption}
Assumption \ref{rip} serves as a widely used identifiable assumption in the literature on high-dimensional statistics \cite{huang2018constructive, Zhu202014241, li2022minimax}. 
It requires that all sub-matrices of $ \boldsymbol{X}^\top  \boldsymbol{X}/n$ with a size smaller than $3s$ exhibit spectral bounds within $[1-\delta_{3s}, 1+\delta_{3s}]$.

\begin{theorem}{($\ell_2$ error bound)}\label{th:convergence1}
 Assume that Assumption  \ref{subgaussian} and \ref{rip} hold.
Let $s\geq s^*$ and $\gamma_{3s} \coloneqq 2\delta_{3s}\sqrt{\frac{1+\delta_{3s}}{1-\delta_{3s}}}<1$. 
Then, with probability at least $1-C_1\exp\{-C_2 s \log (p/s)\}$, we have  
 \begin{equation}\label{eq:first}
  \|{{\boldsymbol{\beta}}^t}-{\boldsymbol{\beta}}^*\|_2 \leq \gamma_{3s}^t \|{\boldsymbol{\beta}}^*\|_2 + \frac{ 12 \sigma}{1-\gamma_{3s}}\left\{\frac{s\log(p/s)}{n}\right\}^{1/2}.
\end{equation}
\end{theorem}
Theorem \ref{th:convergence1} demonstrates that under Assumption \ref{rip}, the convergence rate of HTP is governed by the contraction factor $\gamma_{3s}$. This indicates that the estimation error diminishes geometrically toward the ground truth.
Here, $\gamma_{3s}$ restricts an upper bound for the parameter $\delta_{3s}$, limiting the correlation among the columns of $ \boldsymbol{X}$. 
Furthermore, under the condition of $\delta_{3s}$, $\boldsymbol{\Phi}$ acts as the contraction factor, ensuring a reduction in the optimization error throughout the iterations.
\begin{remark}
    A sufficient condition for $\gamma_{3s}<1$ is $\delta_{3s} < 0.3478$, which is weaker than the condition $\delta_{3s}<0.1599$ in \cite{huang2018constructive} and $0.1877$ in \cite{Zhu202014241}.
\end{remark}


 \begin{corollary}\label{cor:optimal}
Assume that all the conditions in Theorem \ref{th:convergence1} hold. Then, we have
  \begin{align*}
   {\rm{pr}}\left\{\|{{\boldsymbol{\beta}}^t}-{\boldsymbol{\beta}}^*\|_2 \leq \frac{12 \sigma}{1-\gamma_{3s}}\left\{\frac{s\log(p/s)}{n}\right\}^{1/2}\right\}\geq 1-C_1\exp\{-C_2 s \log (p/s)\}
  \end{align*}
  for each $t$ satisfies $
  t \geq \log_{1/\gamma_{3s}} \left\{ \frac{(1-\gamma_{3s})\|{\boldsymbol{\beta}}^*\|_2}{\sigma}\left \{\frac{n}{s\log(p/s)} \right\}^{1/2} \right\}$.
 \end{corollary}
  Theorem \ref{th:convergence1} suggests that the estimation error of HTP is influenced by the non-vanishing terms. Following a few iterations, Corollary \ref{cor:optimal} provides the non-asymptotic upper bounds for the estimation error of HTP. Specifically, when $s$ is of the same order as $s^*$, the HTP algorithm achieves the minimax optimal rate \cite{raskutti2011minimax}.
\subsection{Adaptation to the unknown sparsity and variance}\label{sec:adaptive}
This section addresses the construction of an adaptive estimator that achieves minimax optimality without prior knowledge of the sparsity $s^*$ and noise level $\sigma^2$. 
In the following development of algorithms, we treat $s$ as a tuning parameter and pursue to determine the optimal model size. 
Given an upper bound $s_{\max}$ of $s$, this entails running Algorithm \ref{alg:iht1} along the sequence  $[s_{\max}]$, followed by employing a model selection criterion for identifying the optimal model size.

We use the information criterion which is considered by \cite{verzelen2012minimax}:
\begin{equation}\label{eq:criterion}
  \mbox{IC}(s) = \log \mathcal{L}_n(\hat {{\boldsymbol{\beta}}}^{(s)}) +K\frac{s}{n}\log\frac{p}{s},
\end{equation}
where $\hat {{\boldsymbol{\beta}}}^{(s)}$ is the estimator of Algorithm \ref{alg:iht1} given model size $s$, and $K$ is an absolute constant.
The estimator minimizing the criterion \eqref{eq:criterion} is the optimal model size.
To get an adaptive minimax optimal estimator, the following assumption is required:
\begin{assumption}\label{samplesize}
     The sample size $n$ is large enough such that for some positive constant $C_0$,
	\begin{equation*}
	n \ge C_0 s_{\max} \log \left(p/{s_{\max}}\right).
	\end{equation*}
\end{assumption}
\cite{verzelen2012minimax} demonstrated the impossibility of simultaneous adaptation to $s^*$ and $\sigma$ in the ultra-high-dimensional setting, i.e., $n \leq C_0 s^* \log \left(p/{s^*}\right)$. 
Hence, Assumption \ref{samplesize} is a necessary condition for the adaptation to the unknown sparsity and noise level. Additionally, it offers guidance on setting an upper bound for the maximum model size $s_{\max}$.
\begin{theorem}{(Minimax adaption)}\label{thm:ic}
  Assume that  Assumption \ref{subgaussian} and \ref{samplesize} hold. Assume that $ \boldsymbol{X}$ satisfies \mbox{RIP}$(3s^*, \delta_{3s^*})$ with $\delta_{3s^*} < 0.3478$.
  We run Algorithm \ref{alg:iht1} along sequence $[s_{\max}]$. Denote $\hat{s}$ as the model size identified by information criterion \eqref{eq:criterion}, and $\hat{{\boldsymbol{\beta}}}^{(\hat{s})}$ as the corresponding estimator given model size $s = \hat{s}$.
  Then, we have 
  \begin{equation}\label{eq:adaptive2}
    {\rm{pr}}\left(\hat{s} \leq 2s^*\right)\geq 1-C_2\exp\{-C_3 s^* \log (p/s^*)\},
   \end{equation}
   and
  \begin{align}\label{eq:adaptive1}
    {\rm{pr}}\left\{\|\hat{{\boldsymbol{\beta}}}^{(\hat{s})}-{\boldsymbol{\beta}}^*\|_2 \leq C_1\sigma\left\{\frac{s^*\log(p/s^*)}{n}\right\}^{1/2}\right\}\geq 1-C_2\exp\{-C_3 s^* \log (p/s^*)\}.
   \end{align}
 \end{theorem}
 Theorem \ref{thm:ic} stands as a pivotal outcome of our theoretical discoveries. When combined with the information criterion \eqref{eq:criterion}, the HTP algorithm achieves minimax optimality with unknown sparsity and variance.
 Moreover, Theorem \ref{thm:ic} indicates that the selected model size $\hat{s}$ can be controlled within $O(s^*)$.

\begin{remark}
    In Theorem \ref{thm:ic}, we assume that the RIP-type condition holds for parameter $3s^*$ rather than $s_{\max}$, which imposes a weaker assumption than the Lepski-type method \cite{bellec2018SLOPE}. In particular, we first prove \eqref{eq:adaptive2} by contradiction and control the selected model size $\hat s$ within $2s^*$. Then, based on the result of \eqref{eq:adaptive2}, we only need the RIP condition to hold for parameter $3s^*$.
\end{remark}

\begin{remark}
  Criterion \eqref{eq:criterion} can be considered as a variant of the Birg{\'e}-Massart criterion \cite{birge2001gaussian}.
  It is important to note that \cite{verzelen2012minimax} demonstrated that, using criterion \eqref{eq:criterion}, an estimator derived from exhaustive search with $s \in [(n-1)/{4}]$ attains minimax optimality. However, both \cite{birge2001gaussian} and \cite{verzelen2012minimax} provide a theoretical framework to choose $s$, and do not apply it to a polynomial-time algorithm.
  HTP efficiently procures high-quality solutions for high-dimensional sparse regression, which makes criterion \eqref{eq:criterion} computationally tractable in practice. 
\end{remark}

\section{Oracle rate under the beta-min condition}\label{sec:scaled}

In this section, we establish the theoretical guarantees of HTP under the beta-min condition.
Specifically, in Section \ref{sec:oracle}, we show that HTP can achieve the oracle estimation rate under $\min_{i \in S^*}|{\boldsymbol{\beta}}_i^*| \geq C\sigma\{\log (p/s^*)/n\}^{1/2}$.
Moreover, under a slightly stronger condition $\min_{i \in S^*}|{\boldsymbol{\beta}}_i^*| \geq C\sigma(\log p/n)^{1/2}$, we demonstrate that the HTP estimator converges to the oracle least-squares estimator with high probability.
Section \ref{sec:scaled_ada} investigates the theory of the minimum signal strength of the HTP estimator with different model sizes.
Based on this result, we provide an adaptive tuning procedure to achieve the exact support recovery with unknown sparsity.

\subsection{Oracle estimation rate}\label{sec:oracle}
Consider parameter space $\Omega^p_{s^*, \lambda}$ as 
\begin{equation*}
	\Omega^p_{s^*, \lambda} = \{{\boldsymbol{\beta}}^*\in \mathbb{R}^p : \|{\boldsymbol{\beta}}^*\|_0 \leq s^*\ \text{and}\ \min_{i \in S^*}|{\boldsymbol{\beta}}^*_i| \geq \lambda\}.
\end{equation*}
The beta-min condition $\min_{i \in S^*}|{\boldsymbol{\beta}}^*_i| \geq \lambda$ characterizes the scale of the true signal.
Additionally, we demonstrate that if $\min_{i \in S^*}|{\boldsymbol{\beta}}_i^*| > 2C_1\sigma\{\log (p/s^*)/n\}^{1/2}$, where constant $C_1$ is determined in \eqref{eq:adaptive1}, the selected model size $\hat s$ determined by criterion \eqref{eq:criterion} satisfies
\begin{equation}\label{eq:lam2e0}
    {\rm{pr}}\left(s^*/2 \leq \hat{s} \leq 2s^*\right)\geq 1-C_2\exp\{-C_3 s^* \log (p/s^*)\}.
\end{equation}
 The upper bound of $\hat s$ has been proved in Theorem \ref{thm:ic}. We now prove the lower bound of $\hat s$ by contradiction. Assume $\hat s \le s^*/2$ and denote $\tilde S= S^* \setminus  \text{supp}(\hat{\boldsymbol{\beta}})$, therefore $|\tilde S| \geq s^*/2$ by assumption. Then, we conclude that
        \begin{equation*}
        \|{\hat{\boldsymbol{\beta}}}-{\boldsymbol{\beta}}^*\|_2
        > \|{\hat{\boldsymbol{\beta}}}_{\tilde S}-{\boldsymbol{\beta}}^*_{\tilde S}\|_2
        > C_1 \sigma\left\{ \frac{ s^* \log (p/s^*) }n \right\}^{1/2},
        \end{equation*}
        which leads to a contradiction with \eqref{eq:adaptive1}. Therefore, \eqref{eq:lam2e0} holds.

Result \eqref{eq:lam2e0} asserts that, with high probability, HTP with criterion \eqref{eq:criterion} selects a model with the correct order of dimension under the beta-min condition.
In addition, \eqref{eq:lam2e0} implies that $s^*\in [\hat s/2, 2\hat s]$ with high probability. 
Consequently, 
in what follows we identify the optimal model size within this interval. We only concern the theoretical guarantees of HTP with interval $s \in [s^*,4s^*]$ rather than $[s^*, \infty)$. 

\begin{theorem}{(Oracle estimation rate)}\label{thm:scaled}
  Assume that  Assumptions \ref{subgaussian}, \ref{rip}, and \ref{samplesize} hold.
	Assume that $\delta_{3s} < \epsilon^2 \wedge 0.005$ and $\lambda \geq (1+3\epsilon)\sigma\{2\log (p/s^*)/n\}^{1/2}$ for some $0<\epsilon <1$. 
	Assume that the initial estimator ${\boldsymbol{\beta}}^0$ is a minimax optimal estimator, that is, $\|{\boldsymbol{\beta}}^0-{\boldsymbol{\beta}}^*\|_2 \leq C_1 \sigma \{s^*\log (p/s^*)/n\}^{1/2}$ and $\| {\boldsymbol{\beta}}^0\|_0  \le s$.
	Then, when $s \in [s^*, s^*+s^*/\log (p/s^*)]$, for $t \geq \log \left\{ C_{1,\epsilon}^2  \log (p/s^*)  \right\}$, we have
	\begin{equation}\label{eq:scaled_res1}
  \sup\limits_{{\boldsymbol{\beta}}^* \in \Omega^p_{s^*, \lambda}} {\rm{pr}}\left\{\|{{\boldsymbol{\beta}}^t}-{\boldsymbol{\beta}}^*\|_2 \geq C_{2,\epsilon} \sigma\left(\frac{s^*}{n}\right)^{1/2}\right\} \leq \exp(-C_{3,\epsilon} s^*)+\Big(\frac{s^*}{p}\Big)^{C_{4,\epsilon}}.
	\end{equation}
    Moreover, when $s = s^*+\tau $ where $\tau \in (s^*/\log (p/s^*), 3s^*]$, for $t \geq \log \left( C_{1}^2 s^*/\tau \right)$, we have
	\begin{equation}\label{eq:scaled_res2}
 \sup\limits_{{\boldsymbol{\beta}}^* \in \Omega^p_{s^*, \lambda}} {\rm{pr}}\left\{\|{{\boldsymbol{\beta}}^t}-{\boldsymbol{\beta}}^*\|_2 \geq C_{2,\epsilon} \sigma\left(\frac{\tau \log (p/s^*)}{n}\right)^{1/2}\right\} 
  \leq \exp(-C_{3,\epsilon} s^*)+\Big(\frac{s^*}{p}\Big)^{C_{4,\epsilon}},
	\end{equation}
	where $C_{1,\epsilon}, C_{2,\epsilon}, C_{3,\epsilon}$ and $ C_{4,\epsilon}$ are some positive constants related to $\epsilon$. 
\end{theorem}

\begin{figure}[htbp]
  \centering
    \begin{tikzpicture}
      \draw[<->](9.5,0)--(0,0)--(0,4);
      \draw[red,domain=0.7:3] plot(\x, 0.8) node at (1.8, 1.2){$\frac{\sigma^2 s}{n}$};
      \draw[blue,domain=3:8] plot(\x,{0.5*\x-0.7}) node at (6.8,1.6){$\frac{\sigma^2s \log(p/s)}{n}$};      
      \node[below] at (9.5,0) {$s$};
          \draw [blue, thick, decorate, 
    decoration = {calligraphic brace, 
      raise=3pt, 
      aspect=0.5, 
      amplitude=4pt 
    }] (0,1.8) --  (0,3.3);
    \node[blue] at (-1.5,2.6) {Minimax rate};
    \node[red] at (-1.2,0.8) {Oracle rate};
      \node[below] at (3,0) { $s^*+\frac{s^*}{\log(p/s^*)}$};
      \node[below] at (0.7,0) {$s^*$};
      \node[below] at (8,0) {$4s^*$};
     \node[below] at (5.2,0) {$(1+c)s^*$};
      \draw[dashed](0.7, 0)--(0.7, 0.8);
      \draw[dashed](3, 0)--(3, 0.8);
      \draw[dashed](5, 0)--(5, 1.8);
      \draw[dashed](8, 0)--(8, 3.3);
     \draw[dashed](0, 3.3)--(8, 3.3);
     \draw[dashed](0, 0.8)--(1, 0.8);
     \draw[dashed](0, 1.8)--(5, 1.8);
    \end{tikzpicture}
    \caption{\label{fig:rate}The estimation error bound of HTP estimator with metric $\|\cdot\|_2^2$ for model size $s \in [s^*, 4s^*]$. Here $c$ is some positive constant.}
\end{figure}
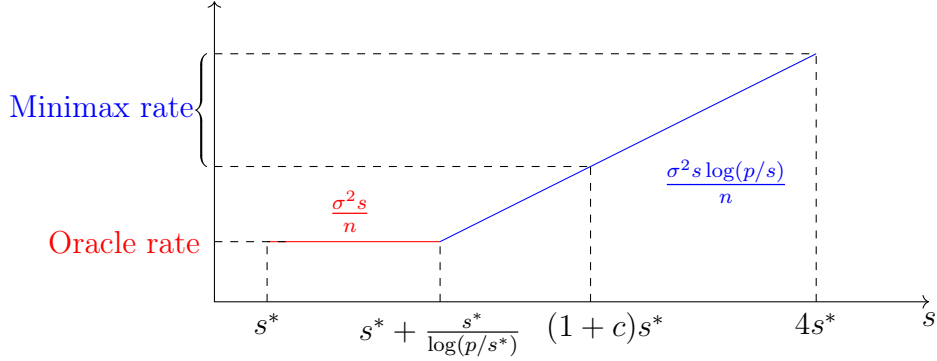

Theorem \ref{thm:scaled} demonstrates that when $s \in [s^*, s^*+s^*/\log (p/s^*)]$, HTP estimator attains the oracle estimation rate $O(\sigma(s^*/n)^{1/2})$ after a few iterations with high probability.
For $s \in (s^*+s^*/\log (p/s^*),  4s^*]$, the estimation error of HTP gradually diverges from the oracle rate and increases to the minimax optimal rate. Figure \ref{fig:rate} illustrates the theoretical results of Theorem \ref{thm:scaled}.

\begin{remark}
\cite{ndaoud2019interplay} established
the oracle estimation rate as 
    \begin{equation}\label{eq:oracle_rate}
\inf_{\hat{\boldsymbol{\beta}}}\sup_{{\boldsymbol{\beta}}^* \in \Omega^{p}_{s^*, \lambda}} {\mathbf{E}}_{\hat {\boldsymbol{\beta}} } \|\hat{\boldsymbol{\beta}}-{\boldsymbol{\beta}}^*\|_2^2 \asymp \frac{\sigma^2s^*}{n},\qquad \forall \lambda \geq (1+\epsilon)\sigma\{2\log (p/s^*)/n\}^{1/2}.
\end{equation}
where $\mathbf E_{\hat{\boldsymbol{\beta}}}$ represents the expectation with respect to $\hat {\boldsymbol{\beta}}$.
As a result, Theorem \ref{thm:oracle} and \eqref{eq:oracle_rate} together show that HTP achieves the minimax optimality under parameter space $\Omega^p_{s^*,\lambda}$.
\end{remark}

\begin{corollary}\label{cor:recovery}
	Assume that all the conditions in Theorem \ref{thm:scaled} hold. 
	Then, when $s \in [s^*, s^*+s^*/\log (p/s^*)]$, for $t \geq \log \left\{ C_{1,\epsilon}^2  \log (p/s^*)  \right\}$, the type-I and type-II errors of variable selection satisfy 
 	\begin{equation*}
\sup\limits_{{\boldsymbol{\beta}}^* \in \Omega^p_{s^*, \lambda}} {\rm{pr}}\left\{\frac{|S^t \Delta S^*|}{s^*} = \Omega\left(\frac{1}{\log (p/s^*)}\right) \right\} \leq \exp(-C_{2,\epsilon} s^*)+\Big(\frac{s^*}{p}\Big)^{C_{3,\epsilon}}.
	\end{equation*}  
\end{corollary}
Corollary \ref{cor:recovery} shows that, when the estimator achieves the oracle rate, both the type-I and type-II errors of variable selection can be controlled within $o(s^*)$ with high probability.

\begin{figure}[htbp]
  \centering
  \includegraphics[scale = 0.8]{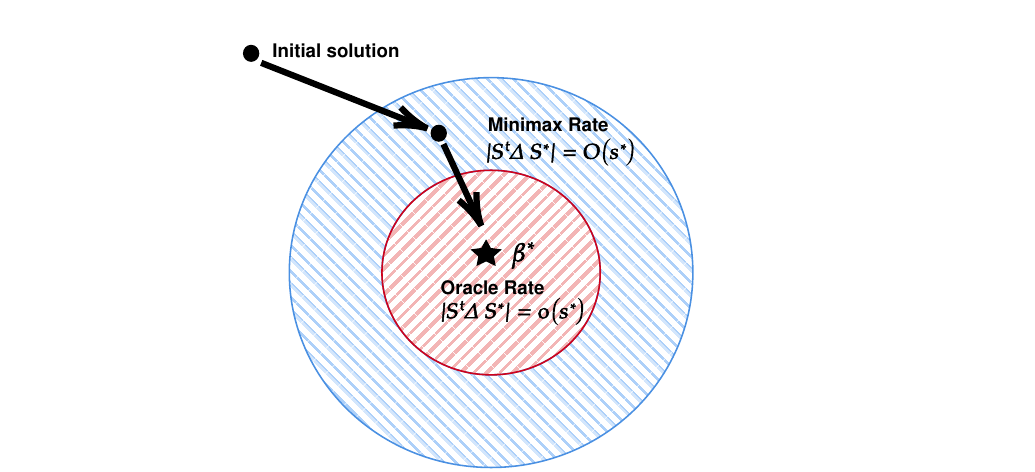}
    \caption{\label{fig:fig4}The solution path of HTP as iteration $t$ increases. The blue region indicates the iterations achieving the minimax rate, while the red region represents the iterations reaching the oracle rate.}
\end{figure}

\begin{remark}
    Figure \ref{fig:fig4} shows the iterative thresholding dynamics of the estimator of HTP. Given an initial solution and letting $s = s^*$, the IHT estimator first hits the minimax optimal region i.e., the blue region in Figure \ref{fig:fig4}. Then, as $t$ increases, the estimator achieves the oracle rate,  i.e., the red region, and the error of variable selection can be controlled at a low level as Corollary \ref{cor:recovery} shows. This behavior benefits from dynamic signal detection, as shown by \cite{butucea2018variable}. Under the beta-min condition, the support recovery result enables HTP to exceed the minimax rate and attain the oracle rate.
\end{remark}

Theorem \ref{thm:scaled} and Corollary \ref{cor:recovery} illustrate that
the HTP estimator performs as well as the oracle estimator on estimation asymptotically, although it may not necessarily recover the true support set $S^*$ exactly. By strengthening the beta-min condition in Theorem \ref{thm:scaled}, Theorem \ref{thm:oracle} demonstrates that the HTP estimator converges to the oracle estimator with high probability. Consequently, Corollary \ref{cor:recovery} is strengthened as Corollary \ref{cor:exact}, affirming that the HTP estimator can exactly recover the true support set $S^*$ with high probability.
Denote the oracle least-squares estimator $\tilde{\boldsymbol{\beta}}^*$ as 
\begin{equation*}
	\tilde{\boldsymbol{\beta}}^*_{S^*} = ( \boldsymbol{X}^\top_{S^*}  \boldsymbol{X}_{S^*})^{-1} \boldsymbol{X}_{S^*}^\top  \boldsymbol{y},\ \text{and}\ \tilde{\boldsymbol{\beta}}^*_{(S^*)^c}=0.
\end{equation*}
In the following theorem, we establish the error bounds between the HTP estimator and the oracle estimator.
\begin{theorem}\label{thm:oracle}
	Assume that all the conditions in Theorem \ref{thm:scaled} hold and $\lambda \geq (2+ 2\epsilon)\sigma(2\log p/n)^{1/2}$.
	Then, given $s = s^*+\tau$ where $\tau \in [0, 3s^*]$, we have 
	\begin{equation*}
\sup\limits_{{\boldsymbol{\beta}}^* \in \Omega^p_{s^*, \lambda}} {\rm{pr}}\left\{\|{{\boldsymbol{\beta}}^t}-\tilde{\boldsymbol{\beta}}^*\|_2 \geq   \left( \frac{9}{10} \right)^t  C_{1,\epsilon} \sigma \left (\frac{s^* \log p  }{n}\right)^{1/2}
		+ C_{2,\epsilon} \sigma\left(\frac{2 \tau \log p}{n}\right)^{1/2}
				\right\} \leq\Big(\frac{s^*}{p}\Big)^{C_{3,\epsilon}}. 
	\end{equation*}
\end{theorem}
For $s > s^*$, Theorem \ref{thm:oracle} demonstrates that the estimation error rate is $O(\sigma((s-s^*)\log p/n)^{1/2})$. Notably, when $s = s^*$, the error diminishes to zero with iterations, implying that the HTP estimator converges to the oracle estimator with a probability approaching 1 as $s^*/p \to 0$.
\begin{corollary}{(Exact support recovery)}\label{cor:exact}
	Assume that all the conditions in Theorem \ref{thm:oracle} hold.  Then, when $s = s^*$ and $t \geq  2 \log \left( (C_2/\epsilon)^2 s^* \right) $, we have
	\begin{equation*}
		\lim\limits_{s^*/p \rightarrow 0} \sup\limits_{{\boldsymbol{\beta}}^* \in \Omega^p_{s^*, \lambda}} {\rm{pr}}\left(S^t = S^* \right) = 1.
    \end{equation*} 
\end{corollary}
Corollary \ref{cor:exact} shows that given $s = s^*$, HTP algorithm can recover the true support set $S^*$ exactly in finite iterations with high probability.
\begin{remark}
   Existing exact support recovery results for $\ell_1$ estimators \cite{zhao2006model, wainwright2009sharp} relied heavily on the irrepresentable condition, recognized as a stronger condition compared to the RIP \cite{van2009}. 
    Typically, achieving exact support recovery demands a tight $\ell_\infty$ error bound, and in this context, incoherence conditions become inevitable. Some studies circumvent these requirements by leveraging the RIP condition, establishing surrogate $\ell_2$ error bounds instead. However, this approach often imposes a more stringent beta-min condition \cite{zhang2018, huang2018constructive,Zhu202014241}.
In our work, we establish the $\ell_\infty$ error bound by asymptotically connecting the oracle ordinary least-squares estimator. Specifically, Theorem \ref{thm:oracle} shows that the HTP estimator asymptotically converges to 
 the oracle estimator when $s = s^*$. Consequently, we can establish a tight $\ell_\infty$ error bound by solely using the RIP condition. Therefore, HTP attains the exact support recovery with the optimal beta-min condition and RIP condition.
\end{remark}

\subsection{Adaptation to the exact support recovery}\label{sec:scaled_ada}
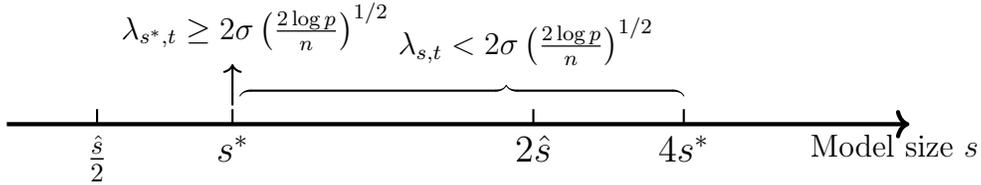
\begin{figure}[htbp]
  \centering
        \begin{tikzpicture}
            \draw[->,ultra thick] (-5,0)--(7,0);
            \node[below] at (6.8,0) {$\text{Model size } s$};
            \node[below] at (-3.8,0) {\large $\frac{\hat{s}}{2}$};
            \draw[black, thick] (-3.8,0) -- (-3.8,0.2);
            \draw[->,thick] (-2,0.25)--(-2,0.8);
            \node[below] at (-2,0) {\large $s^*$};
            \draw[black, thick] (-2,0) -- (-2,0.2);
            \node[below] at (2,0) {\large $2\hat{s}$};
            \draw[black, thick] (2,0) -- (2,0.2);
            \node[below] at (4,0) {\large $4s^*$};
            \draw[black, thick] (4,0) -- (4,0.2);
            \node[] at (-1.7,1.3) {$\lambda_{s^*,t}\geq 2 \sigma\left(\frac{2\log p}{n}\right)^{1/2}$};
            \node[] at (1.9,1.1) {$\lambda_{s,t}< 2 \sigma\left(\frac{2\log p}{n}\right)^{1/2}$};
            \draw [thick, decorate, 
            decoration = {calligraphic brace, 
              raise=5pt, 
              aspect=0.6, 
              amplitude=4pt 
            }] (-1.9,0.2) --  (4, 0.2);
        \end{tikzpicture}
    \caption{\label{fig:lambda}The minimum signal strength of estimator ${\boldsymbol{\beta}}^t$ for different model sizes under the conditions of Theorem \ref{thm:lambda2}.}
    \end{figure}

In this section, we consider an adaptive tuning procedure to identify the optimal model size.
Denote $\lambda_{s, t}$ as the smallest nonzero value of $|{\boldsymbol{\beta}}^t|$ given model size $s$, that is, the minimum magnitude of the nonzero estimate of signal at step $t$ in Algorithm \ref{alg:iht1}.
In the following theorem, we show that after a few iterations, the minimum signal strength of the estimator ${\boldsymbol{\beta}}^t$ exhibits a separable phenomenon for specific model sizes.
\begin{theorem}\label{thm:lambda2}
    Assume that the conditions in Theorem \ref{thm:scaled} hold and $\lambda \geq \big( C_1\vee 4 \big) \sigma(2\log p/n)^{1/2}$, where constant $C_1$ is determined in \eqref{eq:adaptive1}. Given $s = s^*$,  for $t \geq 2\log (C_2^2 s^*)$, we have
    \begin{equation}\label{eq:lam2e1}
        \sup\limits_{{\boldsymbol{\beta}}^* \in \Omega^p_{s^*, \lambda}} {\rm{pr}}\left\{\lambda_{s^*, t} < \big( (C_1-2)\vee 2 \big) \sigma\left(\frac{2\log p}{n}\right)^{1/2}\right\} \leq \Big(\frac{s^*}{p}\Big)^{C_3}. 
    \end{equation}
    On the other hand, given $s \in (s^*, 4s^*]$, for $t \geq 2\log (C_2^2s^*)$, we have
    \begin{equation}\label{eq:lam2e2}
\sup\limits_{{\boldsymbol{\beta}}^* \in \Omega^p_{s^*, \lambda}} {\rm{pr}}\left\{\lambda_{s, t} \geq 2\sigma\left(\frac{2\log p}{n}\right)^{1/2}\right\} \leq \Big(\frac{s^*}{p}\Big)^{C_3}. 
	\end{equation}
\end{theorem}
Theorem \ref{thm:lambda2} demonstrates that when the size $s$ matches the true sparsity $s^*$, the minimum signal strength becomes adequately substantial after $t \geq 2\log (C_2^2 s^*)$. Additionally, for cases where $s > s^*$, the minimum signal strength can be effectively controlled below the threshold of $2\sigma(2\log p/n)^{1/2}$. This relationship is visually depicted in Figure \ref{fig:lambda}.

    To determine the optimal size $s^*$, we suggest an adaptive tuning procedure. A straightforward approach is to apply hard thresholding with a threshold of $2\sigma(2\log p/n)^{1/2}$ to the minimum signal strength of estimator ${\boldsymbol{\beta}}^t$ with model size belonging to $[\hat{s}/2, 2\hat{s}]$. The largest model size for which the minimum signal strength of the estimator ${\boldsymbol{\beta}}^t$ exceeds the threshold $2\sigma(2\log p/n)^{1/2}$ indicates the optimal model size $s^*$.
We call this the adaptive tuning procedure. It can successfully identify the optimal model size $s^*$ with high probability. Combining with Theorem \ref{thm:oracle} and Corollary \ref{cor:exact}, we prove the selection consistency of the HTP estimator.
The HTP estimator adaptively converges to the oracle estimator with high probability.

\begin{remark}
    The optimal threshold of the adaptive tuning procedure relies on the true noise level $\sigma$.
    However, the noise level $\sigma$ is usually unknown in practice. 
    Fortunately, we can obtain a reliable estimate using the minimax optimal estimator $\hat{{\boldsymbol{\beta}}}^{(\hat{s})}$ shown in Theorem \ref{thm:ic}, even in high-dimensional settings. Specifically, we can replace $\sigma$ in the threshold with the plug-in estimator $\hat{\sigma} = \| \boldsymbol{y}- \boldsymbol{X} \hat{{\boldsymbol{\beta}}}^{(\hat{s})} \|_2/n^{1/2}$ to facilitate optimal model size selection.
\end{remark}
\begin{remark}
To enhance its empirical performance, we propose an approach that involves comparing the ratios of the minimum signal strengths of adjacent estimators. By employing a pre-determined threshold $\kappa$, we identify the optimal size by determining the maximum index where the ratio exceeds $\kappa$. In our numerical experiments, we set $\kappa$ to 2.
\end{remark}
In particular, the beta-min condition is necessary to achieve the separation of the minimum signal strength of HTP estimators as shown in Theorem \ref{thm:lambda2}. However, in practical applications, verifying the beta-min condition might be challenging. Hence, evaluating the effectiveness of the adaptive tuning procedure becomes crucial.
Let $\tilde{s}$ be the selected sparsity by the adaptive tuning procedure, and $\hat s $ be the selected sparsity in Section \ref{sec:adaptive}. To assess its efficiency, we examine whether the following condition holds:
\begin{equation*}
\left\|\hat{{\boldsymbol{\beta}}}^{(\tilde{s})}- \hat{{\boldsymbol{\beta}}}^{(\hat{s})}\right\|_2^2 \le C  \sigma^2 \frac{\hat{s}\log(p/\hat{s}) }{n}.
\end{equation*}
If the above condition is met, we conclude that the adaptive tuning procedure is efficient. 
This means that our adaptive tuning procedure is signal-adaptive: when the beta-min condition in Theorem \ref{thm:lambda2} is not satisfied, the estimation error of $\hat {\boldsymbol{\beta}}^{(\tilde s)}$ has the minimax rate $\sigma \{ s^* \log (p/s^*)/n\}^{1/2}$. On the other hand, when the beta-min condition is fulfilled, the estimation error adaptively meets the oracle rate $\sigma ( s^*/n )^{1/2}$. 
We summarize the above observation as the following corollary:

\begin{corollary}\label{cor:efficient}
Given two estimators $\hat{{\boldsymbol{\beta}}}^{(\tilde{s})}$ and $\hat{{\boldsymbol{\beta}}}^{(\hat{s})}$ where $\hat s = \arg\min_{s \in [s_{\max}] } \text{IC}(s)$, and $\tilde s$ is the selected model size by the adaptive tuning procedure. Define an estimator
\begin{equation*}
\hat{\boldsymbol{\beta}} = \begin{cases}
\hat{{\boldsymbol{\beta}}}^{(\tilde{s})}  & \text{, if }\left\|\hat{{\boldsymbol{\beta}}}^{(\tilde{s})}- \hat{{\boldsymbol{\beta}}}^{(\hat{s})}\right\|_2^2 \le C  \sigma^2  \frac{\hat{s}\log(p/\hat{s}) }{n},\\
\hat{{\boldsymbol{\beta}}}^{(\hat{s})}    & \text{, otherwise, }   
\end{cases}
\end{equation*}
where $C = 2 (5+C_1)^2$ and $C_1$ is determined in \eqref{eq:adaptive1}. Assume that all conditions in Theorem \ref{thm:lambda2} hold. Then, as $s^*, p/s^* \to \infty$, with probability tending to 1, the estimator $\hat {\boldsymbol{\beta}}$ achieves adaptivity to the minimum signal strength, that is,
\begin{equation*}
\| \hat{\boldsymbol{\beta}} - {\boldsymbol{\beta}}^* \|_2 \le \begin{cases}
5\sigma \left(\frac{s^*}{n}\right)^{1/2} & \text{, if } {\boldsymbol{\beta}}^* \in \Omega^p_{s^*, \lambda},\\
\left((2C)^{1/2} + C_1 \right)\sigma
\left\{\frac{s^*\log(p/s^*)}{n}\right\}^{1/2}  & \text{, otherwise,}   
\end{cases}
\end{equation*}
where $\lambda = \big( C_1\vee 4 \big) \sigma \left(\frac{2\log p}{n}\right)^{1/2}$.
\end{corollary}
\begin{algorithm}[h]
    \caption{\label{alg:iht2}\textbf {F}ull-\textbf{A}daptive \textbf{HTP} (FAHTP) algorithm}
      \begin{algorithmic}[1]
        \REQUIRE $\boldsymbol{X}$, $\boldsymbol{y}$, $s_{\max}$ and $\kappa$.
        \FOR {$s=1,\ldots, s_{\max}$,}
        \STATE $\hat {\boldsymbol{\beta}}^{(s)}$ = HTP($ \boldsymbol{X},  \boldsymbol{y}, s$).
        \STATE $\lambda_{\min}(s) = \min_{i \in \text{supp}(\hat {\boldsymbol{\beta}}^{(s)})} |\hat {\boldsymbol{\beta}}^{(s)}_i|$.
        \STATE Compute the corresponding information criterion \eqref{eq:criterion} as $\text{IC}_s$.
        \ENDFOR
        \STATE $\hat{s} = \mathop{\mathrm{argmin}}_{s \in [s_{\max}]} \{ \text{IC}_s \}$ and compute $\hat \sigma = \| \boldsymbol{y}- \boldsymbol{X} \hat{{\boldsymbol{\beta}}}^{(\hat{s})} \|_2/n^{1/2}.$
        \FOR {$\tilde s=2\hat{s},\ldots,\hat{s}/2$}
        \IF  {$\lambda_{\min}(\tilde s)/\lambda_{\min}(\tilde s+1) \geq \kappa$ \text{and} $\|\hat{\boldsymbol{\beta}}^{(\hat s)}-\hat{\boldsymbol{\beta}}^{(\tilde s)}\|_2^2 \leq 5 \hat\sigma^2 \hat{s}\log( p/\hat{s})/n$,}
        \STATE Break and let $\hat s = \tilde s$.
        \ENDIF
        \ENDFOR
        \ENSURE $\hat {\boldsymbol{\beta}} = \hat{\boldsymbol{\beta}}^{\hat s}$.
      \end{algorithmic}
    \end{algorithm}

 Consequently, we obtain the asymptotic property of the adaptive tuning HTP estimator as the following corollary.
\begin{corollary}{(Asymptotic properties)}\label{cor:asymp}
    Assume that the conditions in Corollary \ref{cor:efficient} hold. When ${\boldsymbol{\beta}}^* \in \Omega^p_{s^*, \lambda}$, the estimator $\hat {\boldsymbol{\beta}}$ defined in Corollary \ref{cor:efficient} satisfies
\begin{equation*}
		\lim\limits_{ n, s^* \rightarrow \infty,\ s^*/p \rightarrow 0 } {\rm{pr}}\bigg(\hat{\boldsymbol{\beta}}=\tilde{\boldsymbol{\beta}}^*\bigg) = 1.
\end{equation*}
In addition, assume that $s^* = o(n^{1/3})$, $E |\boldsymbol{\xi}_i|^3 \le C_1 \sigma^3$ and $|\boldsymbol{X}_{ij}| \le C_2$ for all $i\in [n],~ j \in S^*$.
Then, for all $\boldsymbol{a} \in \mathbb R^{s^*}$, the asymptotic distribution is
   \begin{equation*}
   \sqrt{n}\cdot \boldsymbol{a}^\top (\hat {\boldsymbol{\beta}}_{S^*}-{\boldsymbol{\beta}}^*_{S^*}) \rightarrow \mathcal{N}\left(0,\sigma^2 \boldsymbol{a}^\top\hat{\boldsymbol{\Sigma}}_{oracle}^{-1} \boldsymbol{a} \right),
   \end{equation*}
    where $\hat{\boldsymbol{\Sigma}}_{oracle} = \boldsymbol{X}_{S^*}^\top \boldsymbol{X}_{S^*}/n $.
\end{corollary}
Therefore, the estimator $\hat {\boldsymbol{\beta}}$ in Corollary \ref{cor:efficient} possesses the strong property when the beta-min condition is met.
Here we set constant $C=5$ and provide the full-adaptive version of the HTP algorithm in Algorithm \ref{alg:iht2}.

\section{Analysis of the general signals}\label{sec:strong_signal}
The previous theoretical guarantees are derived under the assumption of the beta-min condition. However, the beta-min condition may not hold in practice. Therefore, analyzing the error bounds for general signals, i.e., without assuming the beta-min condition, is of particular interest. \cite{fan2023} established a tighter error bound exceeding the minimax rate through implicit regularization with an early stopping criterion. They divided the signals into strong part and weak part, demonstrating that achieving a tighter bound hinges on accurately identifying the number of strong signals and controlling the total magnitude of the weak signals. 

In this section, we divide the support set $S^*$ into $S_1^*\coloneqq\{i \in S^*: |{\boldsymbol{\beta}}^*_i| \geq C_1 \sigma(\log p/n)^{1/2}\}$ and $S_2^*\coloneqq\{i \in S^*: |{\boldsymbol{\beta}}^*_i| < C_1 \sigma(\log p/n)^{1/2}\}$, which correspond to the sets of strong and weak signals, respectively. 
We can derive a generalized version of Theorem \ref{thm:scaled} as follows:
\begin{theorem}\label{thm:signal}
      Assume that the conditions in Theorem \ref{thm:scaled} hold.
	Then, for the above general signal, given $s = |S_1^*|+\tau$, for $t \geq \log \left( C_2^2  \log p  \right)$, we have
 \begin{equation}\label{eq:general1}
 {\rm{pr}}\left\{\|{\boldsymbol{\beta}}^t-{\boldsymbol{\beta}}^*\|_2^2 \geq C_3 \left(\sigma^2\frac{|S_1^*|}{n}+\|\boldsymbol{\beta}^*_{S_2^*}\|_2^2\right)\right\} \leq \exp(-C_{4} |S_1^*|)+\Big(\frac{|S_1^*|}{p}\Big)^{C_{5}}. 
\end{equation}
\end{theorem}
For the general signals, Theorem \ref{thm:signal} characterizes a tighter error bound than the minimax rate.
It demonstrates that the upper bound can be decomposed into two components: the error associated with estimating the strong signals and the error related to the weak signals. Notably, the first component is independent of the dimension $p$, while the second can be controlled by the total magnitude of the weak signals.
A similar result has been achieved through implicit regularization with an early stopping criterion \cite{fan2023}.
Specifically, when $s^* = |S_1^*|$ , we recover the result \eqref{eq:scaled_res1} in Theorem \ref{thm:scaled}. 
\begin{corollary}\label{cor:general_ada}
    Assume the conditions in Theorem \ref{thm:signal} hold and $(|S_1^*|\log p)/{n} \to 0, ~ {|S_1^*|}/{p} \to 0$ and $|S_1^*| \to \infty$. Then, by Algorithm \ref{alg:iht2}, with probability tending to 1, we have
     \begin{equation*}\label{adaptive}
\|\hat{\boldsymbol{\beta}} - \boldsymbol{\beta}^*\|^2_2 \leq \begin{cases}C_2\left(\sigma^2\frac{|S_1^*|}{n} + \|\boldsymbol{\beta}^*_{S_2^*}\|^2_2\right) & \|\boldsymbol{\beta}^*_{S_2^*}\|_2 \leq C_1 \sigma(\log p/n)^{1/2} \\ C_2\left(\sigma^2\frac{|S_1^*|\log p}{n} + \|\boldsymbol{\beta}^*_{S_2^*}\|^2_2\right) & \|\boldsymbol{\beta}^*_{S_2^*}\|_2 > C_1 \sigma(\log p/n)^{1/2} .\end{cases}
\end{equation*}
\end{corollary}
Corollary \ref{cor:general_ada} shows that our proposed FAHTP can achieve a tighter estimation error bound than the minimax rate $\sigma^2s^*\log p/n$ with high probability. It is worth noting that the Corollary \ref{cor:general_ada} can cover the case of $\ell_q$-sparisty ($0 < q < 1)$, which implies that FAHTP is minimax adaptive to the unknown $q$ and radius $R_q$ of the $\ell_q$-ball \cite{raskutti2011minimax}.

\section{Numerical experiments}\label{sec:simulation}
This section provides numerical experiments complementing our theoretical findings and offering insights into the empirical performance of our proposed methods. We begin by demonstrating the convergence phenomenon of the minimum signal strength in Section \ref{num:convergence}. 
Following this, Section \ref{num:adjust} highlights the efficacy of the adaptive tuning procedure under assumptions on the minimum signal strength. Lastly, we conduct a comparative analysis between our proposed FAHTP algorithm and several state-of-the-art methods by synthetic and real datasets.
Here we fix $K = 3$ in the information criterion \eqref{eq:criterion}.

Given estimator $(\hat {\boldsymbol{\beta}}, \hat S)$,
we consider five metrics to measure the efficacy. The first one is the estimation error (EE): $\|\hat {\boldsymbol{\beta}} - {\boldsymbol{\beta}}^*\|_2$. We measure the difference between $|\hat S|$ and true sparsity $s^*$ by sparsity error (SE): $|\hat{S}|-|S^*|$. To evaluate the performance of variable selection, we consider the true positive rate (TPR), false positive rate (FPR) and Mathew's correlation coefficient (MCC). Here MCC is defined as 
\begin{align*}
    \text{MCC} = \frac{\text{TP}\times \text{TN}-\text{FP}\times \text{FN}}{\{(\text{TP+FP)(TP+FN)(TN+FP)(TN+FN)}\}^{1/2}},
  \end{align*}
  where TP= $|\hat S \cap S^*|$ and TN= $|\hat S^c \cap (S^*)^c|$ stand for true positives/negatives, respectively. FP= $|\hat S \cap (S^*)^c|$ and FN= $|\hat S^c \cap S^*|$ stand for false positives/negatives, respectively.

\subsection{The convergence of the minimum signal strength}\label{num:convergence}
We set $n=300$, $p=1000$, $s^*=10$, and $\sigma=1$. Each entry of $\boldsymbol{X} \in \mathbb{R}^{n \times p}$ is independently drawn from $\mathcal{N}(0,1)$.
The non-zero entries of ${\boldsymbol{\beta}}^*$ are independently drawn from $\text{Unif}\left\{2\sigma(2\log p/n)^{1/2}, 10\sigma(2\log p/n)^{1/2}\right\}$.
Subsequently, we run Algorithm \ref{alg:iht1} with model sizes ranging from 1 to 30.
To investigate the convergence phenomenon of the minimum signal strength for each model size, we maintain a fixed iteration number of 20.

\begin{figure}[htbp]
  \centering
  \includegraphics[scale = 0.42]{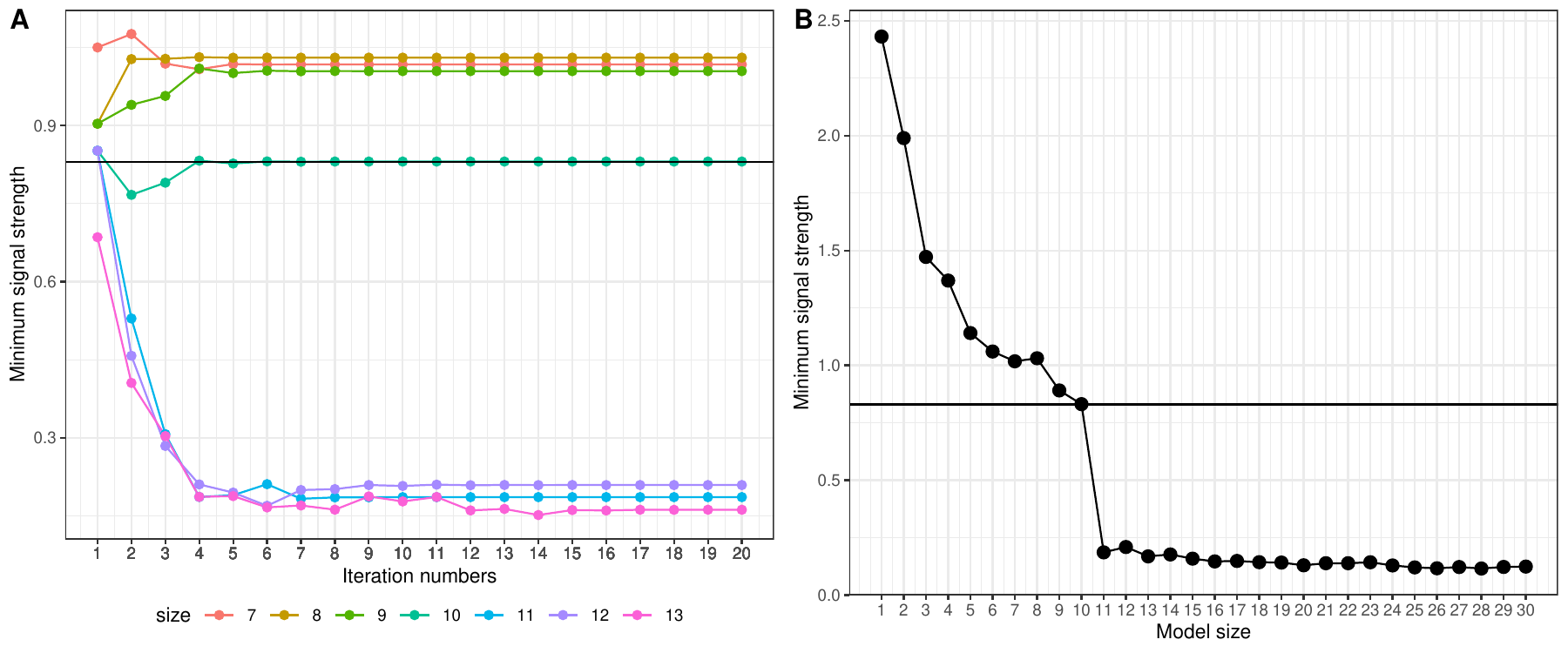}
    \caption{\label{fig:fig1}The phenomenon of the minimum signal strength of HTP estimators with varying model sizes. The horizontal black line represents the true minimal signal strength.
    (A) The convergence of the minimal strength of estimators within 20 iterations.
    (B) The minimum signal strength of estimators with different model sizes.}
\end{figure}
Plot A in Figure \ref{fig:fig1} reveals that after several iterations, all strengths converge. Particularly, when $s=s^*$, the minimum signal strength fluctuates around the true value and eventually converges to it.
Moreover, for $s<s^*$, the minimum signal strength surpasses the true value, whereas for $s>s^*$, it notably falls below the true value. Plot B comprehensively illustrates this phenomenon. Notably, with the assumption regarding the minimal signal strength, there exists a significant gap between model sizes $s=s^*$ and $s>s^*$, aligning with the findings in Theorem \ref{thm:lambda2}. 
This observation indicates the feasibility of applying our adaptive tuning procedure in practice.

\subsection{The benefit of the adaptive tuning procedure}\label{num:adjust}
\begin{figure}[htbp]
    \centering
    \includegraphics[scale = 0.65]{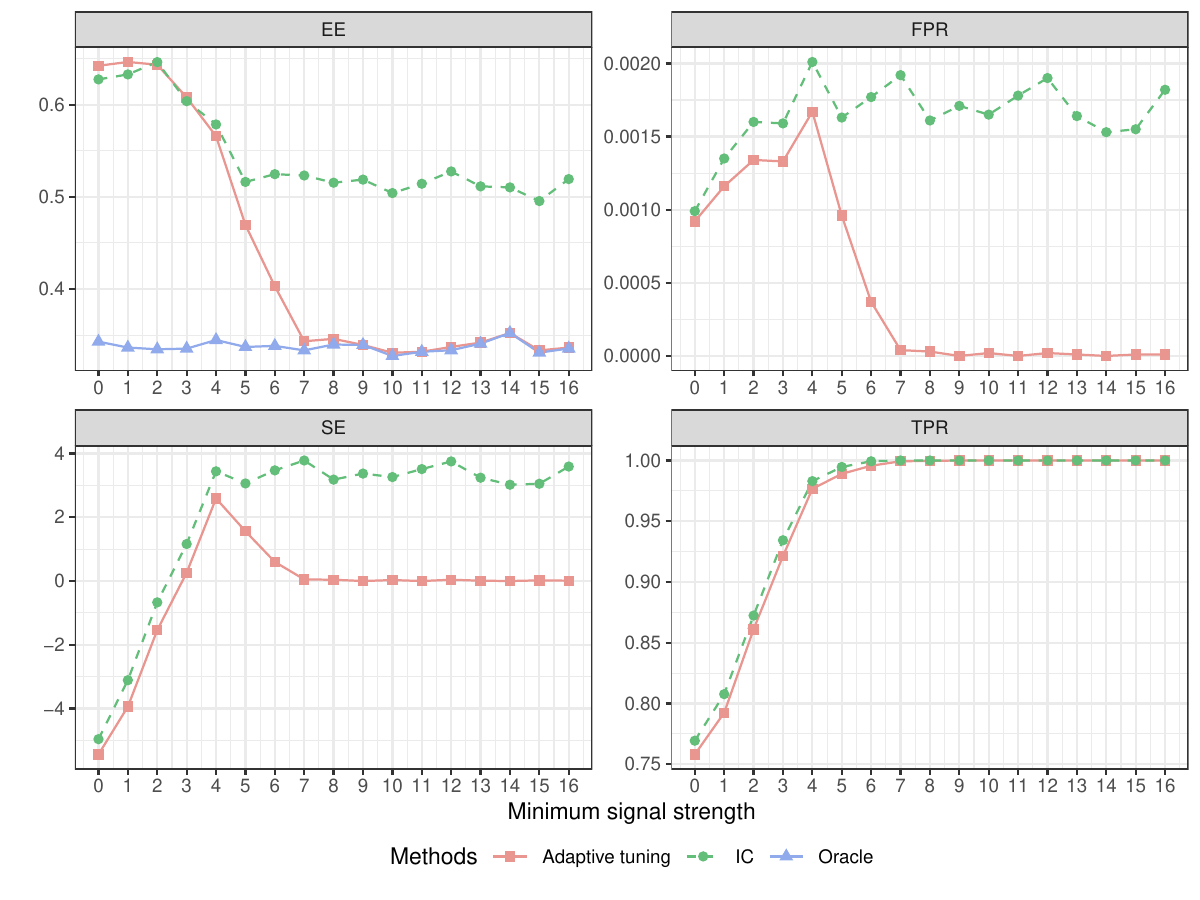}
    \caption{\label{fig:fig2} 
    Performance with increasing minimal signal strength.
    }
  \end{figure}

In this section, we demonstrate the efficacy of the adaptive tuning procedure through simulations.
We set $n=300,~p = 2000,~s^*= 30~,\sigma=1$, and each entry of $\boldsymbol{X}\in \mathbb{R}^{n\times p}$ is independently drawn from $\mathcal{N}(0,1)$. 
To control the minimum strength of the true signal, the nonzero entries of ${\boldsymbol{\beta}}^*$ are independently drawn from $\text{Unif}\left\{ k/4\cdot \sigma(2\log p/n)^{1/2}, 4\sigma(2\log p/n)^{1/2}\right\}$, where $k$ increases from 1 to 16.
The parameter $k$ controls the minimal signal strength, with a larger $k$ resulting in a stronger minimum signal strength.
Here we consider three methods to estimate the signals: 
\begin{itemize}
  \item \textbf{The IC method}: Use information criterion \eqref{eq:criterion} to identify $\hat{s}$.
  \item \textbf{The adaptive tuning method}: Run Algorithm \ref{alg:iht2}.
  \item \textbf{The oracle estimator}: The OLS estimator on the true support set $S^*$.
\end{itemize}
The computational results are shown in Figure \ref{fig:fig2}.

  From Figure \ref{fig:fig2}, it's evident that for small minimum signal strengths, i.e., $k \leq 4$, both the IC and adaptive tuning methods perform similarly, showing significant discrepancies from the oracle estimators.
  As the minimum signal strength increases, the adaptive tuning procedure demonstrates its superiority. It exhibits improved performance in both estimation and variable selection compared to the IC method.
  Specifically, the adaptive tuning procedure effectively limits false positives, resulting in excellent variable selection performance. When the minimum signal strength reaches a sufficiently high level, the adaptive tuning procedure facilitates HTP in achieving the exact recovery of the true support set.
  Remarkably, the estimation error of the adaptive tuning method matches that of the oracle estimator, aligning with our theoretical findings.

\subsection{Comparison with the state-of-the-art methods}\label{num:compare}
We conduct a comparative analysis involving our proposed FAHTP algorithm and several methods: LASSO \citep{T1996}, fitted using the R package $\mathtt{glmnet}$ \citep{JSSv033i01}, and SCAD \citep{fan2001variable} and MCP \citep{zhang2010}, computed using the R package $\mathtt{ncvreg}$ \citep{10.1214/10-AOAS388}.
For parameter tuning across all three methods, we utilize 10-fold cross-validation. Additionally, SCAD and MCP use information criterion \eqref{eq:criterion} to select the optimal tuning parameters. However, due to poor empirical performance, we exclude the results of LASSO tuned by \eqref{eq:criterion}.

\begin{figure}[htbp]
  \centering
  \includegraphics[scale = 0.55]{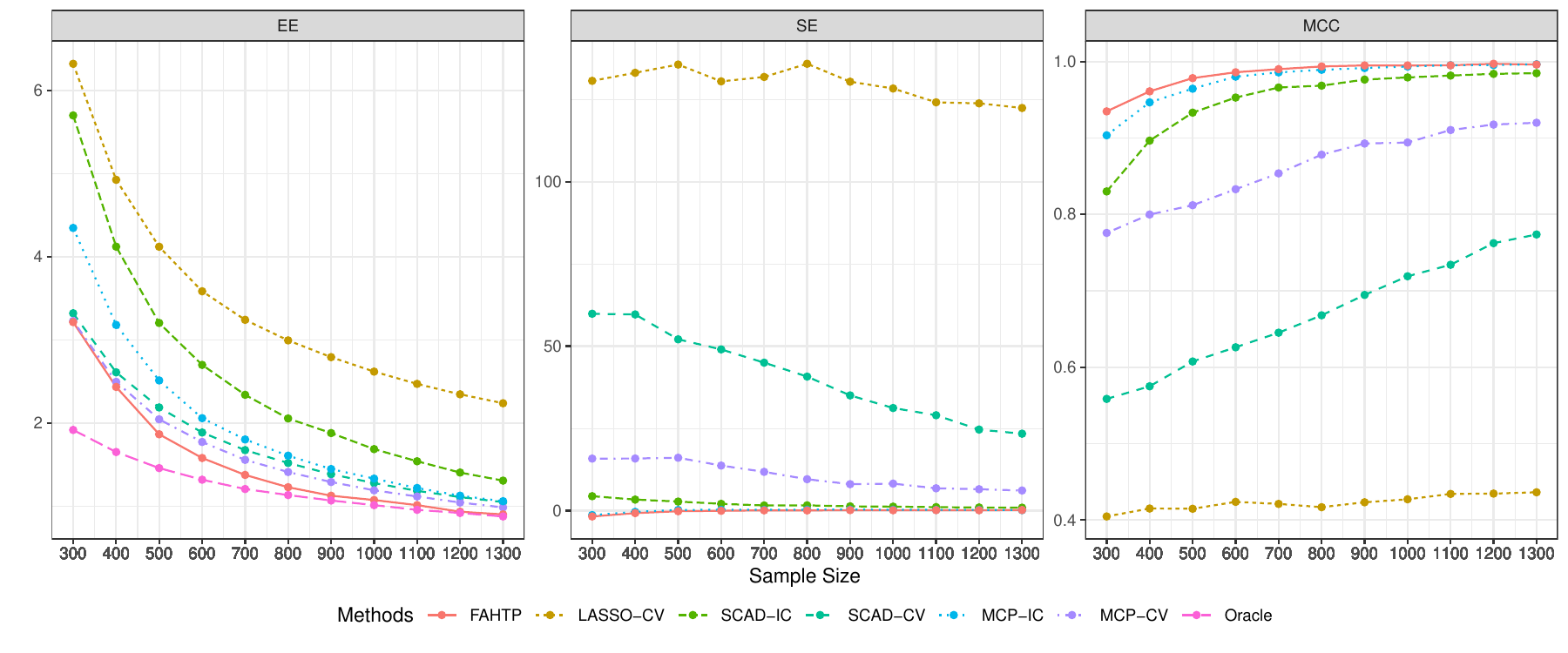}
  \caption{\label{fig:fig3}   Performance metrics with an increasing sample size with 100 repetitions.
  The suffix "-CV" indicates parameter tuning through 10-fold cross-validation. The suffix "-IC" indicates parameter tuning through information criterion \eqref{eq:criterion}.
  }
\end{figure}

We set $p = 2,000, s^*=30$ and increase $n$ from 300 to 1,300 with an increment equal to 100. 
$\boldsymbol{X}$ is generated from the multivariate normal distribution $\mathcal{N}(\mathbf 0, {\boldsymbol{\Sigma}})$, where $ {\boldsymbol{\Sigma}} = 0.5^{|i-j|}$ for $i, j\in [p]$.
${\boldsymbol{\beta}}^*$ are uniformly generated from $[1, 5]$ with a random sign. 
We fix the signal-to-noise ratio (SNR) as 10 to control the noise level $\sigma$, where SNR = $({\boldsymbol{\beta}}^*)^\top {\boldsymbol{\Sigma}} {\boldsymbol{\beta}}^*/\sigma^2$. The computational results are shown in Figure \ref{fig:fig3}.

In terms of estimation error, FAHTP exhibits similar performance to SCAD-CV and MCP-CV when the sample size is small, with all methods displaying substantial deviations from the oracle estimator. However, as the sample size increases, FAHTP notably outperforms the others.
It's worth noting that in this scenario, FAHTP achieves the same estimation error as the oracle estimator. This observation underscores the convergence of FAHTP to the oracle estimator under the assumption of minimum signal strength.
Furthermore, metrics such as SE and MCC demonstrate that FAHTP can exactly recover the true support set with a sufficiently large sample size, reinforcing our theoretical findings.
As discussed in our introduction, while convex methods like LASSO achieve minimax optimality, their estimator bias prevents scaled optimality. Our numerical experiments confirm this phenomenon, showcasing the superiority of FAHTP over convex methods in both theoretical and practical applications.

\subsection{A real data example}\label{num:real}
We illustrate the performance of each method by a supermarket dataset \citep{wang2009}, which consists of 464 records and 6,398 predictors. Each record represents a daily observation from a major supermarket in northern China, with each predictor capturing the sales volume of a specific product. The response variable is the daily customer count. Our goal is to identify which product sales volumes are most strongly correlated with customer numbers.

We conducted 200 random partitions of the dataset. For each partition, 80\% of the records are randomly selected as training data, and the remaining 20\% as the test set. Each method is fitted on the training set, with its performance evaluated on the test set. Table \ref{tab:real} summarizes the average model size and the mean square error (MSE) on the test set. We observe that our method tends to select a sparser model with competitive, even outperformed predictive performance
compared with the state-of-the-art methods. The other methods, based solely on IC or CV, cannot perform well in prediction and variable selection simultaneously.

\begingroup
\renewcommand{\arraystretch}{0.8}
\begin{table}[H]
\centering
\caption{\label{tab:real}The computational results of each method based on 200 random partitions. The standard deviation is shown in parentheses.}
\tabcolsep = 1.5cm
\begin{tabular}{lcc}
\toprule
Method & Model size & MSE \\
\midrule
FAHTP & 2.33(0.53) & \bf{0.410(0.109)}\\
LASSO-CV & 17.51(4.56) & 0.425(0.090) \\
SCAD-IC  & 8.21(1.10) & 0.435(0.119) \\
SCAD-CV & 28.89(15.24) & 0.419(0.109) \\
MCP-IC & 3.90(7.55) & 0.423(0.117) \\
MCP-CV  & 8.05(6.57) & 0.418(0.109) \\
\bottomrule
\end{tabular}
\end{table}
\endgroup

\section{Conclusion and Discussion}\label{sec:conclusion}
HTP is a classical and widely used method for high-dimensional sparse linear regression. However, some theories, such as how to choose the sparsity $s^*$, need to be further studied. This study contributes to the minimax optimality and adaptivity of the HTP algorithm. In this paper, we propose a novel tuning-free method, named FAHTP, to adapt to the unknown sparsity $s^*$. Our procedure ensures the minimax rate $\sigma^2 s^*\log (p/s^*)/n$ first and can achieve the oracle estimation rate $\sigma^2 s^*/n$ under optimal beta-min condition. In other words, we have fully extended the good performance of the hard thresholding estimator in Gaussian sequence models to linear regression, and the assumption on the design matrix solely depends on the Restricted Isometric Property. More importantly, even without the beta-min conditions, our FAHTP achieves a tighter error bound than the minimax rate with high probability.

Notably, this work, especially the conclusion for oracle estimation rate and support recovery, reveals an intriguing aspect of HTP when compared to convex estimators like LASSO, which has been proven that under any signal conditions, it is impossible to break through the minimax optimal rate. Therefore, we theoretically demonstrate some numerical phenomena observed in previous studies where the HTP algorithm outperforms LASSO. These theoretical results inspire us to develop the sparsity parameter selection process.

Moreover, unlike the techniques based on complexity theory and empirical processes in convex algorithms such as LASSO \cite{bellec2018SLOPE}, as well as other non-convex empirical risk-minimizing algorithms such as SCAD and MCP.
\cite{loh2013regularized,fan2014}, our approach combines the analysis of iterative equations with a thorough understanding of the intrinsic properties of the hard thresholding estimator. We think this is the key aspect of our study that enables us to surpass the minimax optimality. To showcase the superiority of the nonconvex algorithm over its convex counterparts, we think that our approach merits significant attention and can be applied to other analyses, such as SCAD and MCP. In summary, the analysis of iterative algorithms can extend the properties of sequence models to practical methods, thereby yielding richer, more detailed, and significant conclusions on statistics. 

\bibliographystyle{unsrtnat}
\bibliography{ref}

\end{document}